\documentclass{article}
\usepackage[utf8]{inputenc}
\usepackage{amsmath}
\usepackage{amssymb}
\usepackage{enumitem}
\usepackage{graphicx}
\usepackage{multicol}
\begin{document}

\title{The Historical Development of Algebraic Geometry\thanks{Special thanks to University of Wisconsin–Milwaukee Department of Mathematical Sciences for approving the content of the video for use in this article.}}
\author{Jean Dieudonn\'e}
\date{March 3, 1972\thanks{This is the original lecture date given with the video. Work on the document began in January 2017, and this version was generated \today.}}
\maketitle

\section{The lecture}
This is an enriched transcription of footage posted by the University of Wisconsin–Milwaukee Department of Mathematical Sciences \cite{youtube}. The images are composites of screenshots from the footage manipulated with Python and the Python OpenCV library \cite{opencvgithub}. These materials were prepared by Ryan C. Schwiebert with the goal of capturing and restoring some of the material.

The lecture presents the content of Dieudonn\'e's article of the same title \cite{dieudonne1972historical} as well as the content of earlier lecture notes \cite{dieudonne1985history}, but the live lecture is naturally embellished with different details. The text presented here is, for the most part, written as it was spoken. This is why there are some ungrammatical constructions of spoken word, and some linguistic quirks of a French speaker. The transcriber has taken some liberties by abbreviating places where Dieudonn\'e self-corrected. That is, short phrases which turned out to be false starts have been elided, and now only the subsequent word-choices that replaced them appear.

Unfortunately, time has taken its toll on the footage, and it appears that a brief portion at the beginning, as well as the last enumerated period (VII. Sheaves and schemes), have been lost. (Fortunately, we can still read about them in \cite{dieudonne1972historical}!) The audio and video quality has contributed to several transcription errors. Finally, a lack of thorough knowledge of the history of algebraic geometry has also contributed some errors. Contact with suggestions for improvement of the materials is welcome.

\subsection{Introduction}
[...] succession of years or periods 
which I have written down on the right. 
I think this perhaps will help 
making things a little clearer. 
I will emphasize in the various periods 
how the various themes come into play 
in the history of the subject.

\subsection{Prehistory}

The first period I which have written there is prehistory 
about 400 BC beginning of geometry 
with the Greeks to 1630 AD. 
Prehistory because in that time 
there is really no algebraic geometry proper.
There is, of course, geometry 
since that is precisely the great invention of the Greeks 
to have invented geometry with proofs, 
which nobody else had ever done before.

They had not only geometry they had also algebra, 
and contrary to some popular beliefs, 
which date back to the 19th century, 
and were rather stupid,
they never tried to divorce algebra from geometry 
as some people think. 
Pure geometry never existed at least: that's an invention of the 19th century.

The Greeks, on the contrary, 
always wanted to do both geometry and algebra at the same time.
They used a lot of algebra to prove theorems of geometry.
As much algebra as they had 
because they hadn't much algebra 
and that was one of their weaknesses.
They didn't evolve a notation 
which would allow them to do algebra 
as was done later on.

Anyway.
They did algebra all right 
to extent they could do it 
at the second degree.
And they used it, 
for instance for the theory of conics 
later on by Apollonius.
And they also used geometry 
to solve algebraic problems

Actually, the invention of the conics, 
contrary to popular belief again, 
which think they were invented 
by cutting a cone by a plane.
This is completely wrong.
The invention of conics came 
through the solution of an algebraic problem, 
namely extracting the cube root of [...].

More precisely they wanted to solve- 
they wrote things as proportion between lengths.
And they wanted to solve 
what they called the double proportion
$a/x=x/y=y/b$.
If you do this via simple algebra, 
you find at once it means finding 
the cube root of $a^2b$.

Well what they did was this: they wrote the equation that we will write now $x^2=ay$ and $xy=ab$.
And Menechmus, a student of Eudoxus 
around middle of the 4th century B.C.,
said ``this is a curve'' and invented the parabola, 
``this is another curve'' and invented the hyperbola,
And cutting both curves gives you solutions. 
And it's only later that they discovered 
this has had to do something with sections of the cone.

So you see, you will say 
the Greeks had some analytic geometry
in the sense that they were able [...]
Certainly they had.
They knew very well how to define a curve 
by a relation between coordinates. 
Except that they did not have algebra 
developed enough to do things 
which were done in later periods, 
in particular to algebraic geometry.

So, they didn't do that only for conics,
they invented a lot of other curves 
both algebraic and transcendental.
As I said, they didn't have any algebra 
so they couldn't explain the difference between these
algebraic and transcendental curves,
nor could they do anything about the [...]
That's why it's prehistory.
Due to the fact they had no algebra, 
they couldn't do really algebraic geometry.
But there was a lot of material which they found.

\subsection{Exploration}
So we come to the second period. 
``Exploration'' I call it,
which is really the beginning of the history.
It starts a sharper beginning, sharp beginning about 1630. 
This is a few years before the simultaneous invention 
by Fermat and Descartes independently of analytic geometry.
Now there was a lot of algebra, 
and they used it immediately to great extent.

They used it for problems which were to become 
the first themes out here:
``Classification'' and ``Transformation.'' 
Both themes are really twins 
because you can't classify things 
unless you can transform the things 
of the same kind to one another in a certain way.

And immediately Fermat discovered... 
explicitly says that... 
he expresses the meaning of dimension.
One curve- one equation in the plane represents a curve.
One equation in 3-space represents a surface and so on.
And the degree, of course, 
Fermat classifies curves of degree 2 
and finds that they are the conics.
And later, Newton classifies the curves of degree 3 
up to change of coordinates and central projections.

In the 19th- in the 18th century, 
they go on doing that with surfaces as well.
They start investigating the first simple fact of singularities [...]  
Curves, double points, also inflection points... things like that.

Without going into much details of course,
they also were interested in intersections of curves.
Newton knows very well 
how to eliminate a variable 
between two algebraic equations polynomials [...]
So he does that when he has two curves of degree $n$, 
and he eliminates one of the coordinates 
and gets from the other an equation, 
which he shows (he discovers) has degree at most $mn$.
This is the first trace 
of what is going to become later 
what we now call B\'ezout's theorem.

B\'ezout himself lives 
at the end of that period around 1760.
And using a better method of elimination than Newton,
he is able to show the equation always has degree $mn$.
But this is still very far 
from what we call multiplicity of intersection,
because what we now call the B\'ezout theorem 
is that to each (point) is attached a number 
which may be any integer [...] one
And the sum of these numbers, the multiplicities, 
must be equal to $mn$.
This is the B\'ezout theorem in its modern form.
So they only have the beginning of that here.

\subsection{The golden age of projective geometry}
So this brings us to the third period
which I call ``the golden age of projective geometry.''
Here again is something which begins rather sharply:
1795 with the date of publication Monge's ``Geometry''
And a little later, about 10-15 years later,
Poncelet develops considerably Monge's ideas.

This is a very sharp- new development,
and of fundamental importance, because
at the same time there are two things which occur.
First, the scalars are extended.
Up to now we only had curves in the real plane 
and coordinates are real numbers.
Now they may be complex values.
The space is extended.
Instead of having points at finite distance, 
we now have points at infinity as well.

And the thing is so successful that for 100 years about [...] 
will do algebraic geometry 
except for things which are immersed 
in the complex projective plane P2(C) 
or the 3 dimensional projective- 
the complex projective space [...]
That will be the thing 
where all geometry will be done for about 100 years.

Of course this had been foreshadowed by some geometers:
projective geometry in the 17th century by Desargues,
complex points by people like Euler or [Stirling?] 
in the 18th century.
But the systematic development dates from Poncelet.

And of course the rewards are immediate.
Things which were quite complicated before 
become extremely simple.
Instead of having five or six types of conics, 
nondegenerate, you have only one. 
in the complex projective plane.
Instead of having 72 types types of cubics, 
which Newton had previously classified,  
we have only three.
And so on.
You get a lot of beautiful theorems 
on conics and quadrics systems of such.

And they start studying also 
the curves and surfaces below degree three and four,
getting very beautiful theorems 
such as the nine points of inflection of a plane cubic,
the 27 lines on a cubic surface,
the 28 double tensions on a plane quartic, 
and so on.
Lots of very beautiful theorems. 

They also start doing a lot about transformation, particularly linear transformations.
They introduce, as you know, the concept of duality,
and immediately they introduce a new invariant. 
They start the degree Now it's a class
A class if of course what hap- 
What the degree is for the dual curve?
When you take it when you transform it by duality,
and then you get a new number
which is of course called the degree in general.
The number of tensions from one point to the curve - the class.

They also begin studying correspondences. 
Not merely one-to-one but one-to-many.
And they get a number of interesting results 
in that direction too.
So this is a beautiful age of development 
of projective geometry.
But it is now more of what we will now call 
linear algebra or linear geometry
than what we now call really algebraic geometry.
This starts in earnest in the next period.

\subsection{Riemann and birational geometry}
The next period is really the birth date 
of modern algebraic geometry.
And it is entirely occupied by one man,
one of the greatest who ever lived: Riemann, 
whose genius and whose importance in this matter 
cannot at all be over-estimated in any sense.

I will try to explain to you what Riemann does in a few minutes.
Riemann starts not at all from geometry, 
although he was perfectly aware … 
that there was a relation with 
what he was doing with geometry [by from] a lattice.
His goal is the study of what are called Abelian integrals.

What are Abelian integrals?
Well, they were not invented by Riemann but by [..] in [date]
You take a curve equation $f(x, y)$ 
in ordinary coordinates, not homogeneous, 
in the plane.
Now remember we are in the complex 
so everything is complex.
And the Abelian integral is an integral like that
What does that mean? 
$R$ is a rational function of two variables, 
so a  quotient of two polynomials.
And you are supposed to replace $y$ 
by the function of $x$ 
taken from this equation.
Put it there.
Then you have a []rational function 
which you integrate - integrate how?
Well, along a path.
You are in the complex 
so you have to go from one point $(a,b)$ of the curve
to another point $(x,y)$ 
by following a path in the complex plane.
If you change the path, 
the integral will change, of course []

And this brings in the concept of periods.
Well, Abel had first introduced these integrals 
in great generality 
--in special cases they were known before-- 
and had tried to classify them.
But Abel had worked 1926.
At that time, the Cauchy theory of integrals 
in the complex plane was just beginning.
So Abel didn't know much of it.
And also Abel didn't know projective geometry, obviously.
He never mentions points at infinity,
and this prevented him from doing 
what Riemann will now do 
in his fundamental papers in 1850 to 1854.

Now what's remarkable in Riemann is that, 
as usually with great mathematicians,
when he takes a problem, 
he takes it in a completely different way from [...]
Doesn't care about curves or about [...]
He starts right away, first of all, 
with the concept of Riemann surface.

You see, what I said here is nonsense.
Take $y$ out of uh equation 
putting in there is something you don't do[...]
There is no such thing as a function with several values.
A function must have one value.
So before doing anything else, 
Riemann starts by inventing the Riemann surface,
so that now, to one point of the Riemann surface, 
there is only one value of the function.
Namely if this equation, say, has degree $n$ in $y$, 
the corresponding Riemann surface is a $n$ sheets complex plane,
each completed with one point at infinity. 
Riemann takes care of that which Abel had not done,
and joined together in the usual way at the ramification points.

Now what is a ramification point?
Well, it is a point 
where you cannot use the implicit function theorem.
It is the derivative of the function if you take it from there is zero.

But Puiseux, just before Riemann started working,
-- Puiseux was a Frenchman, a student of Cauchy -- 
had investigated that in great detail,
and shown that in such a point,
the curve, so to speak, splits 
into a finite number of branches.
And at each branch, you can express the coordinates 
in the neighborhood of the point $X-x_0$ 
the power the positive power 
of some complex variable 
and the other coordinate 
as a power series convergent 
in a neighborhood of zero.
Any function which is meromorphic 
around this point on the Riemann surface [..] 
On a Riemann surface there are
[each?] sheet which are permute 
around the ramification point, 
can be also rewritten as a Laurent series
in powers of $t$ beginning maybe 
with a negative power if we have a pole,
positive power if we have a zero.
So much for the situation before Riemann.

So Riemann starts 
without any reference to algebraic functions
by defining a Riemann surface 
and showing that on the Riemann surface, 
properly understood,
you can do everything which Cauchy did 
in the simple plane.
The whole integration theory of Cauchy 
is now at hand after Riemann's dissertation.
And then he starts with the problem 
of Abelian integrals in his big paper,
in at least the first part of the big paper of Riemann, 
on Abelian integrals.

Well he realizes very well that, 
if he tried to do things in the plane,
you won't get to do anything, 
so he does things on the surface.
And he also realizes that on the surface, 
if you want to do something,
you first have to know how the surface behaves.
That is you have to do algebraic topology: 
analysis situs.
As he [says], 
and that is the first time any mathematician has ever done that,
in advanced algebraic topology on a Riemann surface.

And what he does really 
is the topological classification 
of compact orientable surfaces.
And he shows that on such a surface $S$,
you can draw an even number of simple closed curves
such that when you delete the curves from the surface,
the remaining open set becomes simply connected.
And $2g$ is the minimal number for which you can do that.
So it's an invariant attached to a surface 
which he called a genus.
$g$ is the genus.
$2g$ number of curves.
Genus he showed later is related to the number of ramification points
and the number of sheets by the formula $g-1=w/2 +n$.

OK now after he has done that 
-- though that's already an invention 
of the first magnitude -- 
[..] anything like that
He goes on of course this- 
I must stress here the fact that each of these curves 
has two sides on the surface,
and so you can approach usually in $S'$
--$S'$ is simply connected-- 
you can go from one point to another as you wish. 
And you can also retract the curve on the surface to a point.
But you can always approach a curve $Cj$ from one side or another.

Now Riemann, instead of taking, as I said, 
the problem in the classical way,
states the problem in the following most general way possible.
He says: Let us find meromorphic functions on $S'$
--perfectly well-defined since $S'$ is simply connected--
having the property 
that when you approach the curve $Cj$ from both sides,
the function will not take the same limiting value 
but will have different different values on both sides,
and the difference will be constant.
So, exactly what happens for the Abelian integrals. 
That was known, but the way Riemann puts it 
is a completely new way of attacking the problem.

So this is the problem:
find all meromorphic functions on $S'$ 
having that property. 
These are called integrals of the first or second kind.
There are also integrals of the third kind 
which I won't talk about.
And the way he does it --again completely original--
he says, well, we want to get such a function.
Let's look at the real part of the function.
The real part of the function is a harmonic function[...]
[...]Mereomorphic harmonic function
And it must take given values on each side of the curve
So this is just a simple case says Riemann of Dirichlet's problem
And as he assumed Dirichlet's problem is solved he says we can do it.

You know, the way Riemann thought the Dirichlet problem could be solved 
was fundamentally correct but not proved in his lifetime.
It was only proved about fifty years later.
The Dirichlet problem was something quite difficult,
and Riemann thought he had a way of doing it.
The way didn't work very well.
[..] years later put it on a sound basis.

Anyway, adopting this, the problem is solved.
To each given- you can even take the real part of the period zero
Imaginary parts are not [possible conve- take them possibl- arbitrarily]
And you get your integrals of the first and second kind.
And Riemann proves, in particular,
there are exactly $g$ integrals of the first kind.
First kind means there are no poles.
Second kind there may be poles.
So this is the first of his results.

But now he goes beyond-
and again, as I say, 
he puts the problem completely upside down.
People would have thought of starting 
from well-defined meromorphic functions 
on the surface then integrating them.
He does exactly the reverse.
He starts from integrals,
and says ``now what are the functions?'' 
The meromorphic functions?
Well, they are the functions with which 
you arrive on both sides.
As they are defined on the whole of $S$,
they must take the same value.
In other words, they have no periods.
So he defines meromorphic functions 
on the surface as Abelian integrals 
of the second kind with zero periods.
Then he shows that if you take any two of these,
they are always related by an algebraic relation.
So he gets back to the original point of view
after having gone very far from it.
But a considerable enrichment of course.

Once you have two such functions 
linked by an algebraic relation,
Any other algebraic-- meromorphic function 
is a rational function of $s$ and $t$.
This is proved very easily.
Furthermore the choice of $s$ and $t$ is pretty much arbitrary.
If you make another choice, 
you will have another algebraic relation.
And how do you get from one to another?
Well $s$ and $z$ are rational function of $s_1$ $z_1$ and conversely.
This is what is called a birational transformation.

This was not completely new. Newton, for instance, in the plane, 
had considered the transformation $x' = 1/x$  $y'=y/x$.
This is birational because you can reduce $x$ and $y$ in functions of $x'$ $y'$ and conversely.
[...] change of the degree
Nobody had done anything really 
with these transformations before Riemann,
because nobody had realized 
that something was attached to that.
Invariant number was attached to these transformations:
precisely, the genus.
Two curves, birationally equivalent,
--that is, having Riemann surfaces 
which are isomorphic as- bihomomorphic 
[...] I should say from one to another--
correspond to curves which are birationally transformable 
from one to another.

And this is so important that from the days of Riemann 
until very late in the 19th in the 20th century,
nobody will do anything anymore except birational geometry.
That is, only look at invariants.
The old projective geometry is dead at that time,
and people only are interested in birational geometry.
That is, the degree and class are not anymore invariants, 
but the genuses.

Two curves are in the same genus 
are considered to be the same class because they have-
Well, if they are birationally equivalent.
I should say that of course the genus is not enough 
to determine the curve up to isomorphism.
Actually Riemann himself showed that curves of a given genus or rather isomorphism classes of such curves
Depend on parameters continuous parameters which he called the moduli
$3g-3$ of them [...] 
This was the source of a lot of difficult questions which have only been begun to be solve in our time.

OK, so this is only about half 
of Riemann's paper on Abelian integrals.
You can imagine the wealth of ideas and new methods which are in it.
The other half is just as important 
for [inversion] of integrals,
eta functions, and Abelian varieties, and so on,
But I won't speak about it.
It has had as much influence 
as the first part 
on the rest of mathematics for many years.
So much for Riemann who unfortunately dies 
very young at the age of $40$ in $1866$.

\subsection{Development and chaos}
The next period I call development and chaos.
Development because, of course,
after this fantastic wealth of ideas 
brought by a single man,
in a way which was not always perfectly clear,
analysis was obviously needed correction.

Riemann never bothered about singularities.
He always assumed that curves have no singularities.
But curves do have singularities.
And so on.
This led people in the next period 
to try to put Riemann's discoveries 
on a basis which would be understandable and correct.
Unfortunately, they did not have the genius of Riemann,
and instead of trying to keep in mind 
his fantastic synthesis of analysis geometry algebra
--which he obviously had in mind--
they split it up in two or three different directions.
During that period there will be two or three different schools 
of algebraic geometry who barely understand each other.

I will start with the most interesting of these uh schools for us:
the the so-called algebraic approach.
Not because it is the first historically
--to the contrary it is the last:
it starts in $1882$, 
whereas the other schools start much earlier 
after the death of Riemann--
but because it is the one 
which has had most influence 
after the Riemann papers 
on the development of modern algebraic geometry.

It starts as I say in $1882$ with two big papers:
one by Kronecker which I will barely mention 
although it's a very important paper,
because for algebraic geometry 
it is not of such importance as the other one.
Kronecker essentially puts, 
for the first time, 
on a sound algebraic basis, the concept of dimension,
which was extremely intuitive [..]
He showed that when you define a number of algebraic equations
in, say, 3-dimensional space,
the points that satisfy this equation 
will break up into a union 
of so-called irreducible varieties
of various dimensions well- 
this space can be only 2 dimensions and maybe 
[it will be surfaces no it  $3$ Surfaces curves planes
usable ones]
[...] Kronecker showed
This is relatively not too important.

On the contrary the work of Dedekind-Weber, 
which appeared In the same year,
is of fundamental importance 
because it introduces, 
for the first time, 
the algebraic treatment of curves 
in a birational way.
And Dedekind certainly, next to Riemann, 
is probably the most important man 
in the history of algebraic geometry,
because he also tries to do something 
as great mathematicians do, as I say, 
by taking the problem in a completely different manner 
than the others. 
And what Dedekind and Weber want to do is this: 
they want to have all of Riemann's results 
without any analysis.
Without integrals, without complex, the variables nothing at all.
Just pure algebra.

It was very prophetic 
when we see it from our point of view,
because we know what abstract algebra is.
We've been doing it for fifty years,
but at that time it was completely not understood.
People thought this was not geometry anymore.
They did not understand what Dedekind and Weber were after.
In other words, they were far ahead of their time,
and their importance only became to be recognized 40 years later.

Anyway what they do is this:
they say ``We want to have the Riemann surface,
but the Riemann surface as Riemann considers it 
is something which is linked 
with all sorts of topological or analytic intuitions.
We don't want to have anything to do with that.
We want to have an algebraic definition.''

OK now, as I said, 
one of the fundamental ideas of Riemann
was this birational transformation of curves
which to each curve 
--up to birational transformation-- 
attaches rational functions,
which are always the same, 
since they are the functions of the Riemann surface.
So that is a field.
Really you can add, multiply, divide these functions. 
They still are all rational functions.
So we have a field.

And so what Dedekind and Weber start with 
is that field- and what is that field? 
Well, it is just simply algebraically 
an extension of finite degree 
of the field of rational functions 
in one variable over the complex numbers.
Perfectly well-defined object.
For us it's trivial.
For Dedekind it may not have been quite so trivial,
but he was very familiar with that kind of thing,
since in the ten years preceding that,
he had been doing exactly the same thing 
for the theory of algebraic numbers,
which he had contributed [...]

So we want to go from $K$ to the Riemann surface.
Now how do we get to the Riemann surface?
Let me try to recover how Dedekind may have got to that.

Well, if you have an element of $K$-
what is an element of $K$?
It's a rational function, 
so by Puiseux each branch, branch of the points on the Riemann surfaces,
it has a development which starts 
with a certain number $h$ exponent $h$ 
which may be negative if it has a pole,
positive if its a zero, and so on.
This is called, as you know from ordinary complex variable theory, 
the order of the function.

Let's call it vz0 f at the point z0 which may be- 
if it's not the ramification point of course this is just merely 0
If it doesn't vanish or- 1 or - infinity it's the order of the zero pole.
But you can do it with Puiseux 
even at the ramification point. 
That's very important.
So to each function to each element of $K$ you attach an integer.
And this is quite remarkable that- 
is perhaps the first time in the history of mathematics
that something happened which is now quite commonplace.
When you had- until then when you wrote $f(x)$ for a function, $f$ was fixed and $x$ was variable.
Now this was the first time I think 
in the history of mathematics 
that when one writes $f(x)$, $f$ varies and $x$ is fixed!
Namely to each $z0$ is attached 
for each $f$ --variable function-- a number an integer.
So we have a mapping from $K$ to $z$, 
where the functions are the variable elements 
and the integer is attached to each function at each point.

So that's it if we had the Riemann surface.
But remember, Dedekind does not want to have the Riemann surface 
since he wants to define-
However, having made this observation
he realizes at once that these 
--well, he doesn't say so because he didn't have the word--
these functions have properties 
which now recognized as characterizing what we call a valuation.

A valuation is precisely a mapping from a field 
into the integers which has the two properties:
$w$ of the product is the sum of the $w$'s, 
and $w$ of a sum is at least equal to the minimum of the two.
This is now a very common familiar concept.
Usually it is completed by taking by convention $w(0)=+\infty$
But this is probably the first appearance in mathematics in algebra.

And then they take a big step
What is a Riemann surface?
Each to each point of the Riemann surface,
if we had it,
we would attach a valuation.
They reverse completely the definition and they say:
``The Riemann surface is the set of valuations.''
All possible valuations on [...]
And they start studying this set.

Now to do that,
--roughly very very quickly of course, summarized--
they proceed in the following way.
If they have a valuation of $K$
since $K$ contains a field subfield $C(X)$,
restricted to $C(X)$ it is still a valuation.
Now what are the valuations in the first place on $C(X)$?
Now that is very easy.
There are two kinds.
One which is called the valuation at infinity
which for a polynomial is minus the degree,
and the others
--all this is up to a constant factor--
the others correspond to points of the complex plane.
To each point you attach just simply 
the order of the polynomial at that point.
``$0$'' if it doesn't vanish, ``$1$'' if it has a zero first of all there.
These are the only valuations.

So they correspond really to the geometric picture,
if we had it,
to the plane completed with one point at infinity:
exactly what we what we would have started with.
Furthermore, they showed without too much trouble 
that when you have a valuation here 
and want to extend it to the bigger field,
you can only do it only in a finite number of ways.
So you have the Riemann picture,
but now with completely algebraic language
of the plane completed with one point at infinity,
and above each point there is a finite number of points 
on the Riemann surface.
So the picture fits very well with Riemann's geometric intuition.

And then they start studying this set of valuations of $K$.
I won't go into detail here because I have no time.
They could and probably started studying this 
in a manner which Dedekind had done for algebraic numbers 
by considering the theory of ideals.
But they soon probably recognized 
that this didn't take into account one fundamental fact:
the point at infinity.
The ideals only take care of the valuations 
as they call the finite valuations: 
valuations at finite distance. 
Those that are above points of the plane,
not those which are above the points at infinity.
To take care of the points at infinity, 
they introduce a completely new notion: 
the notion of divisors.

What is a divisor?
We now have the set of valuations, 
which is the abstract surface in the sense of Dedekind,
and a divisor attaches to each valuation an integer,
positive or negative,
which is almost everywhere zero.
So, zero except at a finite number of points.
And of course these form a group.
You can add them.
And you can also order that group.
That means, we say when a divisor is positive.
When a divisor is positive, it means all the alpha v are positive.
Now a positive divisor, of course, is something which is rather intuitive.
It's just simply a finite number of points 
on the Riemann surface 
--on the curve if you wish--
to each of which you attach a positive integer: its multiplicity.
So just a system of points with multiplicities,
geometrically speaking.
Algebraically speaking, it is what I say here.

Now Riemann, when he had [proffered] his paper
which I have barely mentioned,
had considered the following problem:
``Let's take a finite number of points on the curve.
Is it possible to have a rational function
(a function of the field)
which has poles at these points,
with, say, simple poles 
or poles with multiplicities 
bounded from above by given numbers?''
Well, in the plane everybody would say ``the answer is yes of course.''
With one point you can do it,
and always we have a meromorphic function, say $1/z$ 
which has a pole at one point.
Now, Riemann discovered it is not possible in general 
on the Riemann surface.
The total number of points with multiplicities must be at least $g+1$.
$g$ is the genus, plus one.
Then you can do it. 
If it's at least $g+1$ you can do it 
with arbitrary positions of the poles.
If it's smaller than $g+1$ you cannot do it 
except if the poles are at very special positions.
This is a problem which Riemann had investigated.
His student Roch had also completed the results of Riemann 
by precising where these special positions.

And now Dedekind and Weber put that in the following form.
They generalize it by considering, 
for a divisor,
all the elements of the field 
which have valuation at least equal to $-\alpha v$
(the choice of minus is for convenience)
for all valuations.
Remember the $\alpha v$ are almost always zero,
but there's some of them may be positive,
so you have something which allows for poles.
And they show, without too much trouble, that this is a vector space,
and it is finite dimensional.
Finite dimension is called small L of d.

There is also an equivalence relation between divisors.
Two divisors are equivalent 
if their difference is what is called a principal divisor.
Principal divisor attaches to each $f$ in the field 
the divisor where you take as positive elements for $v$ 
the points where the function is zero
with the order of the zeros
and negative elements the points 
where the function has poles [...] order of the poles.
It turns out that the equivalence of divisors 
keeps the degree and the value of L of the dimension, 
as it's called L(d), the same.

And Now this problem is solved by Dedekind-Weber, 
the Riemann-Roch theorem, in their formulation:
$L(D)-L(\delta-D)=\deg(D)+1-g$
Delta is a specific class of divisors, 
what is called the canonical class.
It is linked with the differential forms on the Riemann surface, 
and it is perfectly well-defined.

Now let's see, for instance, 
how this gives back the Riemann theorem about the poles.
If we want to have, say, simple poles at $g+1$ points,
this means that we must take $\alpha v$ equal plus one 
at z $g+1$ points and zero elsewhere.
So the degree is $g+1$.
$g+1 +1 - g$  that makes $2$.
This is a positive or zero number so $L(d)$ is at least $2$.
Now also you always have a constant which satisfies this property.
The valuation of a constant is always zero.
So it's certainly bigger than $-1$.
So the constants are there, and if the dimension is $2$,
that means there are functions in that space L(d) 
which are not constants.
That is exactly what Riemann wanted: 
functions which have poles at $g+1$ given points.
Instead if you took $g$ instead of $g+1$ this wouldn't work,
because this would give only $1$ and not $2$.
So much for this algebraic approach of Dedekind and Weber.

Now before that in 18th after Riemann's death,
other German geometers, 
such as Brill, Noether, and Clebsch and Gordan, 
had started doing similarly an interpretation of Riemann's results
but only in geometric terms.
And we may describe the idea of all these people
by saying they didn't want anything to do with negative divisors.
Divisors were not invented,
but they worked as if there were only positive divisors.
That is, systems of points.

Now, how did they get to this equivalence relation 
which is very fundamental?
Well, they showed that if you have two positive divisors [...]
(begins drawing overlapping curves)
Riemann had already shown that in special cases.
On the curve, that is, systems of points,
they may be obtained by cutting the curves 
by auxiliary curves called adjoint,  
which depend linearly in their equations on parameters.
And if you change the parameter, you get another system of points,
which turns out to give the positive divisor 
equivalent to the preceding one.

So we have a geometric picture 
of what are divisors positive and equivalence, 
and it is with this picture that this school of geometers
did what Dedekind did.
It can be done.
It's a bit more complicated.
More geometric but less clear from the algebraic point of view.
Something which they also did

Two things must be mentioned here.
First of all, they studied singularities 
much more than had been done before.
In the birational picture,
you can replace a curve by a birational transform.
Now in their way of doing things,
they wanted to study not only curves without singularities,
as practically Riemann had done,
but also curves with singularities.
So the problem was:
``Is it possible to transform birationally a curve with singularities 
into a curve which has no singularities?''
And they showed that such was indeed the case.

Similarly they introduced a new concept
which until then was not very widespread:
the idea of not working in three dimensional projective space
but in $n$ dimensional projective space.
After this time, everybody works in $n$ dimensional space.
And then it can be shown that any curve 
can be birationally transformed
into a curve in some higher dimensional space
without any singularities.
So the point of view of Riemann was perfectly justified.
One can work with curves without singularities at all.

Next come the extension of all this 
to higher dimensional varieties
and first of surfaces.

Now Clebsch and Max Noether, in 1870,
had already thought of going in the way Riemann had taken. 
That is, instead of considering a curve, 
they considered a surface in 3- complex plane 
in nonhomogenous coordinates,
and double integral this times of rational functions
taken with the same [...]
That is, extract $z$ from this one and carry it here.
You have a bunch of $x$ $y$ and you have to integrate it.
But the situation became rapidly more complicated,
because what does that mean, this integration?
Remember it's a double integral over a 2 dimensional [...] 4
image in the complex plane of some triangle, say, or something like that.
So immediately we are now in 4 dimensional real space 
(2 dimensional complex space),
and we have there to study what happens 
on integrals over 2 dimensional real varieties
or chains or whatever you like.
This obviously cannot be done intuitively anymore.
And even more when you get to 3 dimensional complex varieties
in 6 dimensional real space.
You need something more.

Although Clebsch, Gordan, and a little later Picard studied these such integrals
and also simple integrals of exact differentials
(always complex, remember)
it's only after Poincar\'e in 1895 had started to invent the simplicial machinery
necessary to start doing algebraic topology in arbitrary dimensional space 
that the thing could get moving.
Well, he himself used to great extent his topology theory
--algebraic topology--
42:04
to get results for algebraic surfaces.

Clebsch and Noether had found that 
there is an invariant quite similar to Riemann's genus 
attached to any surface.
Namely, the number of integrals of this pole which are the first kind.
It is everywhere finite,
but when they tried another way of getting to the genus 
by means of the adjoint surfaces,
now they found another number.
So there were now two genus for a surface:
the geometric genus defined this way
and the arithmetic genus, 
which I won't try to define now because 
the classical definition is extremely complicated and difficult to understand.

Anyway, they proved that the difference was always positive,
and they called it the regularity of the surface.
And for a long time they tried to prove the regularity of the surface
was half of the first Betti number of the surface,
in the sense of algebraic topology,
and also the number of simple integrals of the first kind.
It's only Poincar\'e in 1910 after 15 years of efforts
by Picard and the Italian geometers
who succeeded in proving that.

\begin{figure}[ht]
\caption{}
\label{fig:changing5}
\centering
\includegraphics[width=0.80\textwidth]{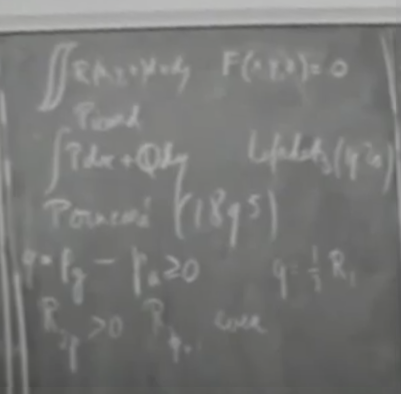}
\end{figure}

And later, Lefschetz continued these ideas of Poincar\'e
and extended them to higher dimensional algebraic varieties
around 1920
and obtained a lot of beautiful results about the Betti numbers.  [See figure \ref{fig:changing5}] 


For instance the even Betti numbers are always different from zero.
The odd Betti numbers are even, if I may say so.
Betti numbers of odd dimension are even integers. 
They may be zero.
And so on.
Plenty of beautiful theorems which Lefschetz found
by using the methods of Poincar\'e using algebraic topology.
It turns out that they are applicable
because Poincar\'e had already sketched a proof,
which was later proved on a very sound basis,
that any algebraic variety is triangularable.
So you can do algebraic topology on it without [any trouble?]

Finally, I should say a word about the Italian school.
The Italian school wanted to do for surfaces and higher dimensional varieties
what Brill and Noether had done for curves.
Namely they wanted to do uniquely geometric constructions and study them.
Of course divisors can be defined just in exactly the same way.
The fact that it's a curve here doesn't matter.
You just simply replace this field 
by a field of rational functions of two variables or three variables.
Otherwise the definitions are the same.

But now instead of having systems of points on the surface,
you have systems of curves,
because a point is something where a function is zero, for instance, at zero.
Well systems of points a function is zero on a surface
That means they would have a system of curves.
On three dimensional variety we have a system of surfaces, and so on.
So you have to do the things which Brill and Noether did 
with systems of points
with systems of curves or surfaces.

Now it turns out that it is far more difficult to work with.
You can manage as Brill and Noether had done
to work with linear systems of points,
but linear systems of curves and surfaces introduce a lot of complexity,
because it can very well be that that divisors
--naturally defined such as this canonical divisor
which can also be defined for surfaces--
are not positive anymore.
This one is a positive divisor for curves. 
For surfaces it may well be not positive.

So the Italians were led to introduce 
all sorts of complex definitions.
Very hard to understand
What they call virtual definitions something like that
And usually well it's almost impossible to read now,
because it goes about that way.
They start giving definition of such and such a thing all right [...]
And then ten lines later they say: 
``Oh yes, but in certain cases you should not define it that way.
You should define it slightly differently 
because there are some troubles there.''
And ten lines later again there is another way.
After two pages it is utterly impossible 
to understand what they have defined.

Well,  what is amazing is that with this kind of handicap
They managed to unearth a large number 
of extremely beautiful and deep theorems
in the theory of surfaces [...]
They got all sorts of invariants.
Not only these two pg pa,
but also what Enriques has called the \emph{plurigenus} $P_k$,
which can have any value integer from $2$ on,
and proved a number of beautiful theorems.
For instance if a surface has arithmetic genus and second plurigenus zero,
then it is rational: 
it is birationally equivalent to the complex plane.
The complex projective plane, I should say.
It's an ordinary theorem due to Enriques.
And similarly all sorts of theorem of that kind.
It's really marvelous 
that people with such handicaps 
were able to do such beautiful things.

\subsection{New structures in algebraic geometry}
And so we come to the sixth period 
Beginning in 1920
What I call new structures in algebraic geometry
Well up to now you have seen 
that algebraic geometry is like the god Janus.
It has two two faces:
on one side it faces complex variable theory with Riemann,
on the other side with Dedekind 
it faces algebraic number theory and 
what we now call commutative algebra.

(47:30)
In that period beginning in 1920,
these analogies become structural paths of algebraic geometry.
They become part of the theory itself and not merely analogies.
And I have listed only two of the --I cannot go into detail--
but I have listed two of these [advance] which happened during that period.

First of all, what I call ``the return to Riemann''
for higher dimensional surfaces
via the idea of K\"ahlerian manifold.
I'll try to explain that in a few words.
When Poincar\'e had started his research on algebraic topology,
he had observed that the exterior differential forms
must have some relation with homology of the surface.
Exterior differential forms are defined, of course,
in local coordinates as given by exterior products of differential coordinates
with coefficients [...]

It's only in 1931 that de Rham was able 
to put this intuition of Poincar\'e 
in a precise form in the famous de Rham theorem
which says that the homology of the surface $Hj$ with 
respect to respect to a manifold we call $M$ 
with coefficients in the real field 
is paired with the cohomology of the exterior form.
Cohomology of the exterior form is defined by looking at the complex form
consisting of the form of degree one degree two et cetera degree $n$
[...]dp
with the differential being the usual exterior derivative.
This gives rise to the cohomology groups,
and they turn out to be paired in a natural duality 
with the homology groups of the manifold.

Now a complex algebraic variety,
in which without singularities, I assume,
is a differentiable manifold, clearly, but much more.
It has a structure of a complex manifold. 
That is, locally you have complex coordinates,
and the transition functions are [holomorphic?]
So what you can do with that, first of all,
is to divide your $p$ forms into finer classes,
what is called forms of type $r$-$s$ with $r+s=p$.
This means simply that instead of taking 
the real and imaginary parts as local coordinates,
you take the complex variables themselves as local coordinates,
and you take the differential of $r$ of the $z$'s and $s$ of the $z$ bar.
These are called forms of type $r$-$s$ [...]

Well, there is even more than that.
There is what is called, on the surface, structure of a K\"ahlerian manifold.
This is a special type of Riemannian manifold, which I won't try to describe here,
but which has a great importance in the study of [the integral?]
The man who recognized that is Hodge in a series of papers starting in 1930's
and which I call ``the return to Riemann''
because what Hodge did was exactly this:
on the Riemannian manifold --you see that Riemann here appears in the second half
it is a difference [...] it is the Riemann of the Riemannian manifolds and not of the Riemann surfaces [...]--
on the Riemannian manifold of arbitrary dimension,
Beltrami had shown that you can define an operator
that is exactly the same thing 
as the generalization of the usual Laplacian in ordinary space.
So you can define harmonic functions.
However, what Hodge did something much less trivial.
He showed you can define harmonic forms on the Riemannian manifold,
--which I will not describe here either,
but this was exactly the generalization of Riemann's idea--
but this time on $n$ dimensions.
Because he proved the fundamental theorem
that in each cohomology class of de Rham, 
there is, on the Riemannian manifold, one and only one harmonic form.

So essentially this rather abstract cohomology group was essentially
isomorphic in a natural way with the space of $j$-harmonic forms,
since there is one and only one in each cohomology class.
So, this already was a natural generalization of Riemann
because it shows, for instance,
that you can always find a harmonic form which has prescribed periods
over independent homology types 
--homologically independent [sites?]--
which is exactly a generalization of Riemann's result,
and proved that the intervention of the harmonic functions in Riemann
was not a [toy?] device or a trick,
it was really deep fundamental principle.
Riemann had guessed it [...]

Now the merit of Hodge is to have done that in the first place.
But that's not the only thing.
It was then he passed to the K\"ahlerian manifold
--the manifold with this special property--
and observed that you can,
in this harmonic form,
in that case they split according to the type.
For each type of have a space of harmonic forms,
so that $Hp$ turns out to be the direct sum of all harmonic forms of type $r$-$s$
$r+s=p$.

Now this has a lot of fundamental consequences.
For instance, it's quite clear that you take a form of type $r$-$s$
and take its complex conjugate.
What happens?
Well the $\bar{z}$ becomes $z$ and the $z$ becomes $\bar{z}$,
so you get a form of type $s$-$r$.
Now this goes on to cohomologies just as well,
and therefore you have the dimensions of $Hrs$ and the dimension of $Hsr$ are the same.
Now suppose $p$ is odd.
Then there are as many $rs$ as there are $r$
So what do you say?
The dimension of that space is even:
that's one of Lefschetz's theorems.
Similarly you can prove all the other Lefschetz theorems
by interpreting them in this theory of Hodge of harmonic integrals.
It is a very beautiful result [to return to].

\begin{figure}[ht]
\caption{}
\label{fig:changing1}
\centering
\includegraphics[width=0.83\textwidth]{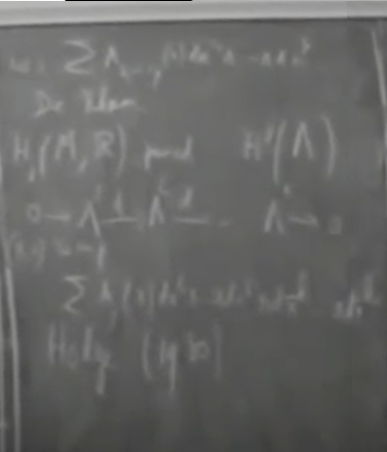}
\end{figure}

Now I come to a completely differently other face of algebraic geometry:
the algebraic face.
And this is not, as we will see, a return to Riemann but a return to Poncelet. [See Figure \ref{fig:changing1}]




You know that in the early 1900s, essentially the German school
had developed what we now call algebra,
which was called abstract algebra at that time.
Namely, that they they realized algebra had nothing to do with real numbers 
or any kind of number.
It was just simply a set of structures 
--groups rings fields et cetera-- of arbitrary elements.
And that was [...] the deep idea behind algebra.

So they they started doing it 
and it and it was immediately apparent
that if you could do all that kind of thing 
all algebra on arbitrary elements,
you could also do algebraic geometry.
Because, after all, what is algebraic geometry?
It is essentially the study of equations 
that form a polynomial in $x_1, x_2, x_n$.
Now you can just simply take the coefficients of the polynomial
and the values of the variable to be in a field $K$,
and then you have an algebraic variety defined over the field $K$.
Now it was quite natural to study these things.
In addition, there were powerful influences which lead to that also.
Namely, coming from algebraic number theory.
But unfortunately it would take an entire lecture for that.

Anyway, around 1920 the German school essentially
started studying algebraic geometry over an arbitrary field.
The people were essentially Emmy Noether (the daughter of Max Noether)
Krull, van der Waerden beginning around 1925,
and later on Zariski and Andr\'e Weil after 1940.

What kind of problem was there?
Well, the essential problem was this.
Up to then, except for Dedekind and Weber,
Riemann the Italians, Brill, Noether everybody even Poincar\'e
had worked with a kind of geometric intuition
founded on essentially on analysis [...] justified.
And the Italians in particular 
had used that with a kind of recklessness
which other people thought was rather extraordinary,
a bit exaggerated.

There were two things in particular that the Italians used
and which were subject to discussion.
One was the concept of points in general positions.
Well, usually when they had a theorem to prove
the geometric algebraic geometric objects which existed
depended on parameters.
And they said ``We will prove the theorem when the points are in general position.''
That means ``when the arguments work.''
And then they said once it's proved for general position, it's always true.

Now this may sound silly,
but a little later Severi,
one of the best representatives of the Italian school,
justified that in 1912 in a very simple way.
He said: ``Your object depends on complex parameters.
These complex parameters usually form a certain set 
in a high dimensional space [...]
And this set has an interior.
Well, general position just means that the parameters are in the interior.
And now everything which we proved depends continuously on parameters,
so when we reach a boundary point, 
if the thing was true in the interior
it's still true at the boundary point.''
This is essentially correct,
although the Italians never bothered to express that in so many words.
They say the thing is in general position.

Alright another thing which was a bit more difficult
was the question of multiplicity of intersection,
which at that time was on a sound basis for curves,
but the Italians used it 
for intersections of surfaces curves and higher dimensional varieties
also rather freely.
For instance, they did not hesitate to speak, on a surface,
of the number of intersections of a curve with itself.
This sounds absurd because a curve with itself has infinitely many points,
so there can't be [...]

Anyway it's correct [...] already shown
by the following argument.
In Poincar\'e's theory of algebraic topology
you can define chains or cycles
and assign them intersection numbers [oriented?].
Well, it's not always defined.
It's not defined a priori when the two are equal for instance.
But when they have complementary dimension 
and have general position, 
in the sense which is perfectly rigorous with Poincar\'e's definition,
you can very well count the points of intersection 
with with value plus or minus one.
It's not always plus one: it may have negative sign.
Take the sum of that,
and the fundamental property of this number 
is that it is invariant when you form $U$ and $V$ 
so as to have the same homology class.
It depends on the homology classes of $U$ and $V$ 
and not $U$ and $V$ themselves.
So this is a perfectly well-defined.

And then you can define it for arbitrary representatives of a cohomology class
In particular, a surface. 
Perfectly meaningful in algebraic topology.
We define the set intersection of a 1-cycle with itself
because the idea you just simply
take a homology which deforms a little bit $C$
and take the intersection in the usual sense,
plus or minus values you have a number.
This is perfectly correct.

Now Severi argued that this was what algebraic geometry did.
Essentially the argument was not very clear,
because he had not made perhaps the observation
which Lefschetz made a little later.
Namely, that instead of taking arbitrary cycles,
you take algebraic varieties.
Then it turns out that this number which may be plus or minus before, 
it was always plus. 
There is a natural orientation for algebraic complex manifolds,
and if you take that orientation, the number of intersections
is always plus one, never minus one.
So the number in the sense of Poincar\'e
was exactly the number of the algebraic geometers' sense,
which had always been considered
and therefore you could very well,
as Severi argued, 
do that for arbitrary algebraic cycles
by taking the topological definition of Poincar\'e as intersection number.
And therefore it was perfectly justified to speak 
of intersection of a curve with itself,
in that sense.
Always positive numbers occur [...] algebraic values

OK, but that was alright as long 
as we were in the complex field,
because there we had the topological theory of Poincar\'e.
What about an arbitrary field?
Nothing like that existed.
So it had to be done something.
And the man who did most to inaugurate this 
was van der Waerden,
who, as I said, got to his ideas 
by returning to ideas of Poncelet.

\begin{figure}[ht]
\caption{What is the intersection number?}
\label{fig:intersection}
\centering
\includegraphics[width=0.80\textwidth]{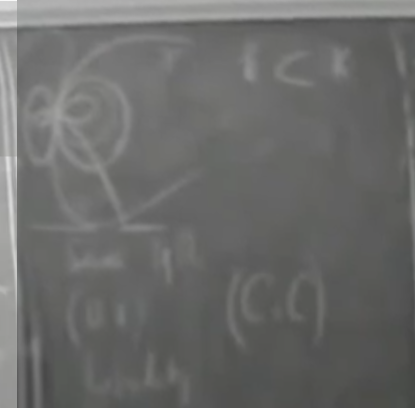}
\end{figure}

Now what had Poncelet done?
Poncelet had had the following idea.
When it came to intersection multiplicities of two curves,
B\'ezout's theorem for instance,
we have two curves which intersect.
Let's say very badly they have 
a bad kind of singular point,
and the other one is just as bad.
It goes that way and touches many times. [See Figure \ref{fig:intersection}]
What is the intersection number?
What is the multiplicity of that intersection?
Nobody knows.

\begin{figure}[ht]
\caption{A path of a continuous family of curves}
\label{fig:curvefamily}
\centering
\includegraphics[width=0.80\textwidth]{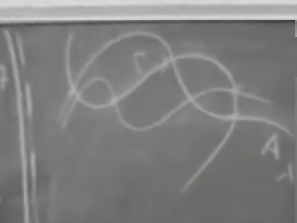}
\end{figure}

Poncelet's idea is this:
let one of the curves fixed (gamma)
and vary the other one.
Well, a curve is always a path of a continuous family of curves,
say, all of the curves of the same degree, for instance.[See Figure \ref{fig:curvefamily}]
[Vary?] the other one
so when it moves
and at certain moment, it will be so nice 
that it will cut the first curve 
in a finite number of simple points with distinct tangents.
In that case of course the number of intersections is $1$.
Now Poncelet says ``now let that parameter go back to its original position
which is so bad.
The number of multiplicity is the number of points which will tend to that point
when the parameter returns [to the path?]''
A very natural definition,
which a good definition.

Now what van der Waerden did
was to generate all points 
he got that-- how did he do that?
For generate point, he went back to Poncelet for some other idea.
Namely, the idea of extending the field.
That had been lost aside.
It had been lost completely,
because people had gone from the reals to the complex field
and been so happy in the complex field for 100 years,
and they never wanted out.
Complex field has all the virtues
It is an algebraically closed field.
It has the theory of complex analytic functions in it.
Everything which makes it the most beautiful field you can work with.

However, here when you started with an arbitrary field $K$,
the idea was to get to a bigger one.
Namely, an equation with coefficients in $K$
still has meaning if the coordinates are in a bigger field $K$.
So you could define,
for a given algebraic variety,
the points of that variety in a bigger field.
And it turns out this was necessary to define what was
to replace the points of general position of the Italians
and generate points.
I won't try to define it here.

Similarly for the intersection multiplicity
Van der Waerden took his ideas from this Poncelet idea, of course,
varying the things.
And he could show that there is a definition of multiplicity
which is perfectly correct.
The only trouble is that with his definition,
multiplicity was not always defined,
even when the intersection where the varieties were intersecting very properly.
Andr\'e Weil a little later was able to remove that
by more fine devices and defined intersection multiplicity 
in every case where it has meaning.
And then it was possible to do as Poincar\'e had done
To define what is the product of two cycles [footage ends]

\section{The chalkboard}
\subsection{THEMES}

\begin{multicols}{2}
\begin{enumerate}[label=\Alph*)]
\item Classification
\item Transformation
\item Infinitely near points
\item Extending the scalars
\end{enumerate}
\columnbreak

Invariants dimension

degree [...] genus

Correspondence

Morphism

Singularities

Multiplicity of intersections

Nilpotent elements

Complex points

generic points

Change of basis

\end{multicols}

\begin{figure}[ht]
\caption{Themes A-D}
\centering
\includegraphics[width=0.83\textwidth]{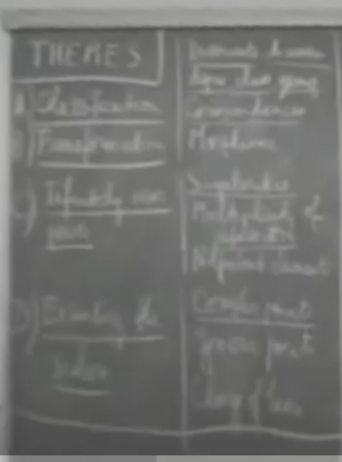}
\end{figure}

\begin{multicols}{2}
\begin{enumerate}[label=\Alph*)]
\setcounter{enumi}{4}
\item Extending the space
\item Analysis and topology in algebraic geometry
\item Commutative algebra and algebraic geometry
\end{enumerate}
\columnbreak

projective geometry

$n$-dimensional geometry

[abstract varieties]

[schemes]

[...]

[...]

K\"ahlerian manifolds

fiber bundles

sheaves

Field of rational functions

Valuations, ideals 

divisors, local rings 

Homological algebra
\end{multicols}

\begin{figure}[ht]
\caption{Themes E-G}
\centering
\includegraphics[width=0.75\textwidth]{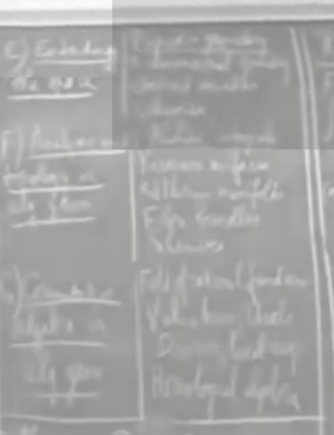}
\end{figure}

\subsection{Central panels}

\begin{figure}[ht]
\caption{Dedekind-Weber panel}
\label{fig:section4}
\centering
\includegraphics[width=0.75\textwidth]{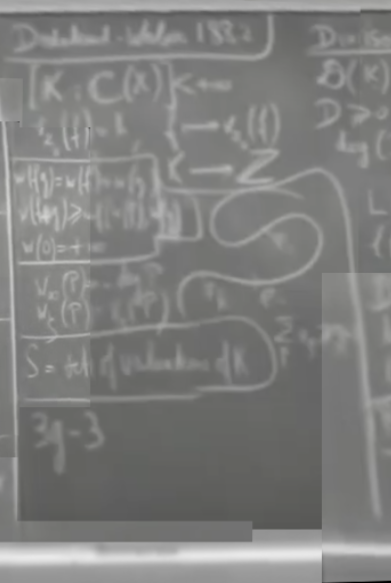}
\end{figure}
See Figure \ref{fig:section4}

$[K: \mathbb C(X)] < \infty$

$r_2(f)...$  $f\to v_a(f)$
$w(fg)=w(f)..(g)$
$v(f+g)\geq ()$
$w(0)=+\infty$
$w_\infty(P)=deg P$
$w_S(P)=..P$
$S=\text{ set of valuations of $K$}$
$3g-3$

\pagebreak
\begin{figure}[ht]
\caption{Dimension panel}
\label{fig:section5}
\centering
\includegraphics[width=0.75\textwidth]{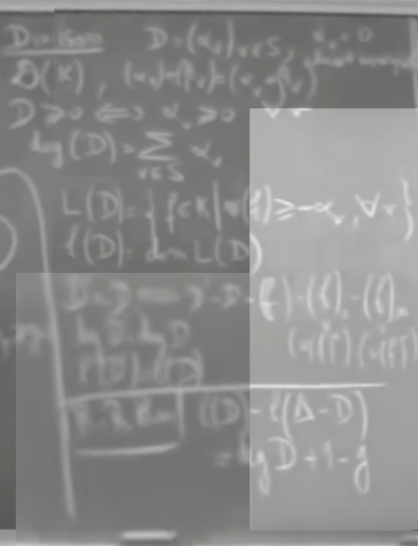}
\end{figure}
See Figure \ref{fig:section5}

Dimension  $D=(\alpha_i)_{i\in S}$, $\alpha_o = 0$ almost everywhere
$\mathcal{D}(K)$ $(w_\alpha\ast(\beta_w)=(w_alpha-\beta_v)$
$D\geq 0 \iff \alpha _v\geq 0$ $\forall v$
$deg(D)=\sum_{v\in S}\alpha_v$
$L(D)=\{f\in K\mid o(f)\geq -\alpha_v, \forall v\}$
$\ell(D)=dim L(D)$
$D'\sim D\iff D'-D=(f)-(f)_a-(f)_b$
by $D$ by $D$ $(v(f)(v(f))$
$f(D) =\ell(D)$
R-R Thm 
$\ell(D)-\ell(\Delta-D)=\deg(D)+1 - g$

\subsection{PERIODS}
\begin{enumerate}[label=\Roman*.]
\item Prehistory (-400-1630)
\item  Exploration (1630-1795)
\item  The golden age of projective geometry (1795-1850)
\item Riemann and birational geometry (1850-1866)
\item Development and chaos (1866-1920)
\begin{enumerate}[label=\alph*)]
\item Algebraic approach 
\item Brill-Noether
\item Integrals on higher dimensional varieties
\item The Italian school
\end{enumerate}

\item New structures in algebraic geometry (1920-1950)
\begin{enumerate}[label=\alph*)]
\item K\"ahlerian manifolds and return to Riemann
\item Abstract algebraic geometry
\end{enumerate}

\item Sheaves and schemes (1950-)
\begin{enumerate}[label=\alph*)]
\item The Riemann Roch theorem and sheaf cohomology
\item Serre [...] varieties
\item Schemes and topologies
\end{enumerate}
\end{enumerate}

\bibliographystyle{IEEEtran}
\bibliography{./citations}

\end{document}